\documentclass[12pt]{article}    

\input{amssym.def}
\input{amssym.tex}

\begin{document}      
\title{On embedding of the Bratteli diagram into a surface}    

\author{Igor  ~Nikolaev\\
Department of Mathematics\\
2500 University Drive N.W.\\    
Calgary T2N 1N4 Canada\\
{\sf nikolaev@math.ucalgary.ca}}

 \maketitle    
    
\newtheorem{thm}{Theorem}    
\newtheorem{lem}{Lemma}    
\newtheorem{dfn}{Definition}    
\newtheorem{rmk}{Remark}    
\newtheorem{cor}{Corollary}    
\newtheorem{prp}{Proposition}    
\newtheorem{exm}{Example}    

\newcommand{\N}{{\Bbb N}}
\newcommand{\F}{{\cal F}}
\newcommand{\R}{{\Bbb R}}
\newcommand{\Z}{{\Bbb Z}}
\newcommand{\C}{{\Bbb C}}   
 \begin{abstract}    
We study $C^*$-algebras ${\cal O}_{\lambda}$ which arise in
dynamics of the interval exchange transformations and measured
foliations on compact surfaces. Using Koebe-Morse coding
of geodesic lines, we establish a bijection between Bratteli
diagrams of such algebras and measured foliations. This
approach allows us to apply $K$-theory of operator
algebras to prove strict ergodicity criterion and
Keane's conjecture for the interval exchange transformations.

\vspace{7mm}    
    

\vspace{5mm}
{\it AMS (MOS) Subj. Class.:  19K14, 46L40, 57R30.}
\end{abstract}

\section*{Introduction} 
Let $\lambda=(\lambda_1,\dots,\lambda_n)$ be a partition of the unit interval
into a disjoint union of open subintervals. Let $\varphi: [0,1]\to [0,1]$
be an interval exchange transformation (with flips). 
Consider a unital $C^*$-algebra ${\cal O}_{\lambda}$ generated by the
unitary operator $u(\zeta)=\zeta\circ\varphi^{-1}$ and characteristic
operators $\chi_{\lambda_1},\dots, \chi_{\lambda_n}$ in the Hilbert
space $L^2([0,1])$.  This (noncommutative) $C^*$-algebra has 
an amazingly rich geometry.

${\cal O}_{\lambda}$ is Morita equivalent to a groupoid $C^*$-algebra  corresponding to 
measured foliations on a compact surface of genus greater than one.
To this end, ${\cal O}_{\lambda}$ is an extension of the irrational rotation
algebra $A_{\theta}$ whose theory experienced an extensive development in the
last decades.

${\cal O}_{\lambda}$ is closely related to simple $C^*$-algebras of
minimal homeomorphisms on a Cantor set. These $C^*$-algebras
were in focus of a brilliant series of works of I.~F.~Putnam  starting
with the papers \cite{Put1}, \cite{Put2}. We refer the reader to
our work \cite{Nik2} for discussion of connections between Putnam's 
algebras and ${\cal O}_{\lambda}$.

The $K$-groups of ${\cal O}_{\lambda}$ are finitely generated 
and can be obtained from the Pimsner-Voiculescu diagram for
the crossed products. Namely, $K_0({\cal O}_{\lambda})=\Z^n$,
$K_1({\cal O}_{\lambda})=\Z$, where $n$ is the number of intervals
in the partition of $[0,1]$. 
The dimension group $(K_0,K_0^+,[1])$ of ${\cal O}_{\lambda}$
was calculated in \cite{Nik2}. (The reader is referred to
Appendix for the details of this construction.) 
When $\varphi$ is minimal the dimension group $(K_0,K_0^+,[1])$
is simple.

Recall that state on dimension group is a positive homomorphism 
of $(K_0,K_0^+,[1])$ to $\R$ which respects the order units $[1]$ and $1\in\R$.
The state space
$S_{\bullet}$ of $(K_0,K_0^+,[1])$ is a Choquet simplex of dimension
$\le n-1$. The dimension of $S_{\bullet}$ is equal to the number of
linearly independent invariant ergodic measures of $\varphi$.
Each invariant measure corresponds to a $1$-dimensional linear subspace
of the state space and $dim~S_{\bullet}=1$ if and only if the interval exchange 
transformation $\varphi$ is strictly (uniquely) ergodic.

It might be one of the most intriguing
problems of topological dynamics since 25 years to indicate
conditions of strict ergodicity of $\varphi$. Some results in this
direction are due to Veech and Boshernitzan.
In 1975 Keane conjectured that ``typical'' $\varphi$ is strictly ergodic.
Masur \cite{Mas} and Veech \cite{Vee}  proved this conjecture in positive using  methods
of complex analysis and topological dynamics, respectively.

This note is an attempt to study dynamics of $\varphi$ using the
ideas and methods of operator algebras. A foundation to such an
approach is given by the following main theorem (to be proved
in Section 2):
\begin{thm}\label{range}
Let $n\ge 2$ be an integer. Let $(P,P_+,[u])$ be a simple
and totally ordered
\footnote{Total ordering condition ensures that the Unimodular Conjecture
is true, see Effros \cite{E} and Elliott \cite{Ell}. 
The author believes the condition is technical, but cannot
drop it at this stage.}
dimension group of order $n\ge 2$. Then there exists an interval
exchange transformation $\varphi=\varphi(\lambda,\pi,\varepsilon)$
of $n$ intervals and a $C^*$-algebra ${\cal O}_{\lambda}$ with the
group $(K_0,K_0^+,[1])$ which is order-isomorphic to $(P,P_+,[u])$.
The transformation $\varphi$ is minimal.
\end{thm}
The proof of the above theorem is based on the identification
of the infinite paths of Bratteli diagram  
with the symbolic  geodesics on a compact surface (so-called
Koebe-Morse theory). This method has an independent interest
since it provides direct links between geometry of geodesics
and $K$-theory of operator algebras.

\bigskip\noindent
The paper is divided into  five sections. In Section 1 we
introduce notation and a lemma on positive cone in 
$K_0({\cal O}_{\lambda})$.
In Section 2 we give the proof of main theorem.
In Sections 3 and 4 we apply Theorem \ref{range}
to establish a strict ergodicity criterion and Keane's
Conjecture, respectively.
Section 5 is an Appendix
containing quick review of dynamics of the interval exchanges,
measured foliations, $K$-theory and rotation numbers associated to the
$C^*$-algebra ${\cal O}_{\lambda}$. The reader is encouraged
to read ``Conclusions and open problems'' section at the end of
this paper.

\medskip\noindent
{\bf Acknowledgements.} I wish to thank G.~A.~Elliott for many helpful 
discussions and ideas. 



\section{Notations}
Let $A$ be a unital $C^*$-algebra and $V(A)$ 
be the union (over $n$) of projections in the $n\times n$ 
matrix $C^*$-algebra with entries in $A$.
Projections $p,q\in V(A)$ are equivalent if there exists a partial
isometry $u$ such that $p=u^*u$ and $q=uu^*$. The equivalence
class of projection $p$ is denoted by $[p]$.

Equivalence classes of orthogonal projections can be made to
a semigroup by putting $[p]+[q]=[p+q]$. The Grothendieck
completion of this semigroup to an abelian group is called
a {\it $K_0$-group of algebra $A$}. 

Functor $A\to K_0(A)$ maps a category of unital
$C^*$-algebras into the category of abelian groups so that
projections in algebra $A$ correspond to a ``positive
cone'' $K_0^+\subset K_0(A)$ and the unit element $1\in A$
corresponds to an ``order unit'' $[1]\in K_0(A)$.
The ordered abelian group $(K_0,K_0^+,[1])$ with an order
unit is called a {\it dimension (Elliott) group} of $C^*$-algebra $A$.

For the $C^*$-algebra ${\cal O}_{\lambda}$ one easily finds
that $K_0({\cal O}_{\lambda})=\Z^n$, see the Appendix. It is
harder to figure out the positive cone $K_0^+({\cal O}_{\lambda})$.
The rest of the section is devoted to this specific question.

\bigskip\noindent
Let us fix the following notation:

\vskip0.5cm
\halign{\indent#\hfil&\quad#\hfil\cr
$\Bbb H$               & Lobachevsky complex half-plane $\{z=x+iy|y>0\}$\cr 
                       & endowed with the hyperbolic metric $ds=|dz|/y$;\cr
                       &\cr
$\partial {\Bbb H}$    & absolute, i.e. line $y=0$ of the Lobachevsky half-plane;\cr 
                       &\cr
$G$                    & Fuchsian group of the first kind;\cr 
                       &\cr
$M_{g,m}$              & orientable surface of genus $g$ with $m$ boundary components;\cr        
                       &\cr                                     
${\cal F}$             & measured foliation of $M_n$ obtained as suspension over\cr
                       & interval exchange transformation $\varphi=\varphi(\lambda,\pi,\varepsilon)$
                         with $n$ intervals;\cr
		       &\cr	 
$\Lambda$              & geodesic lamination corresponding to ${\cal F}$.\cr 
		       &\cr	 
$\gamma$               & geodesic ``generating'' $\Lambda$, i.e. $\bar\gamma=\Lambda$.\cr }

\vskip0.5cm\noindent
Thurston has shown that each measured foliation $\cal F$ can be represented by
a ``geodesic lamination'' $\Lambda$ consisting of disjoint non-periodic geodesics,
which lie in the closure of any one of them; cf Thurston \cite{Thu}. 
Denote by $p: {\Bbb H}\to M_{g,m}$ a covering mapping corresponding to the
action of discrete group $G$.

The geodesic lamination $\Lambda$ is a product
$K\times {\Bbb R}\subset M_{g,m}$, where $K$ is a (linear)
Cantor set. The preimage $p^{-1}(\Lambda)\subset {\Bbb H}$ 
is a collection of geodesic half-circles without self-intersections
except, possibly, at the absolute. The ``footpoints'' of these half-circles
is a subset of $\partial {\Bbb H}$  homeomorphic to $K$.

 Fix a Riemann surface $M_{g,m}={\Bbb H}/G$ of genus $g$ together with a point $p\in M_{g,m}$.
 Let $\gamma$ be a ``generating'' geodesic of the lamination $\Lambda$, i.e. such 
 that closure $\bar\gamma=\Lambda$. Consider the set
\begin{equation}
Sp ~(\gamma)=\{\gamma_0,\gamma_1,\gamma_2,\dots\},
\end{equation}
of periodic geodesics $\gamma_i$ based in $p$, which monotonically 
approximate $\gamma$ in terms of ``length'' and ``direction''. 
The set $Sp ~(\gamma)$  is known as {\it spectrum} of $\gamma$
and is defined uniquely upon $\gamma$.

Let $n=2g+m-1$. Then the (relative) integral homology 
$H_1(M_{g,m}, \partial M_{g,m};{\Bbb Z})\cong {\Bbb Z}^n$. 
Since each $\gamma_i$
is a 1-cycle, there is an injective map $f:Sp~(\gamma)\to 
H_1(M_{g,m}, \partial M_{g,m};{\Bbb Z})$, which relates every closed geodesic its
homology class. Note that $f(\gamma_i)=p_i\in {\Bbb Z}^n$
is ``prime'' in the sense that it is not an integer multiple of some
other point of lattice ${\Bbb Z}^n$. Denote by $Sp_f(\gamma)$ 
the image of $Sp~(\gamma)$ under mapping $f$. Finally, let 
$SL(n,{\Bbb Z})$ be the group of $n\times n$ integral matrices
of determinant $1$ and $SL(n,{\Bbb Z}^+)$ its semigroup
consisting of matrices with strictly positive entries. 
It is not hard to show, that in appropriate basis in 
$H_1(M_{g,m}, \partial M_{g,m};{\Bbb Z})$ the following is true: 

\medskip
(i) the coordinates
of vectors $p_i$ are non-negative;

\smallskip
(ii) there exists a matrix
$A_i\in SL(n,{\Bbb Z}^+)$ such that $p_i=A_i(p_{i-1})$ for
any pair of vectors  $p_{i-1},p_i$ in $Sp_f(\gamma)$.

\begin{dfn}\label{positivecone}
The ordered abelian group
$({\Bbb Z}^n,({\Bbb Z}^n)^+,[1])$ defined as inductive limit
of simplicially ordered groups:
\begin{equation}\label{d.g.}
\Z^n\buildrel\rm A_1\over\longrightarrow \Z^n
\buildrel\rm A_2\over\longrightarrow
\Z^n\buildrel\rm A_3\over\longrightarrow \dots,
\end{equation}
is called associated to geodesic $\gamma$.
\end{dfn}
(We have shown in \cite{Nik2} that the order structure on $({\Bbb Z}^n,({\Bbb Z}^n)^+,[1])$
is independent of the choice of $M_{g,m}$ and $\gamma$.)
\begin{lem}
dimension group $(K_0,K_0^+,[1])$ of the $C^*$-algebra ${\cal O}_{\lambda}$
is order-isomorphic to the associated group $({\Bbb Z}^n,({\Bbb Z}^n)^+,[1])$
of Definition \ref{positivecone}.
\end{lem}
{\it Proof.} See \cite{Nik2}.
$\square$

\section{Proof of Theorem 1}
Let us outline main idea of the proof. To every 
dimension group $(P,P_+,[u])$ with $P\simeq\Z^n$ one can relate a Bratteli diagram
$(V,E)$. The path space $X$ of $(V,E)$ can be made to a topological
space by putting two paths ``close'' if and only if they coincide
at the initial steps. ($X$ is called Bratteli-Cantor compactum.)
$X$ can be embedded (as topological space) into the complex plane 
$\Bbb H$ by identification of each $x\in X$ with a geodesic in $\Bbb H$ 
via Morse coding of the geodesic lines. We show that $X=p^{-1}(\Lambda)$,
where $\Lambda$ is Thurston's geodesic lamination on the surface 
$M_n={\Bbb H}/G$; cf Thurston \cite{Thu}. A concluding step is to recover 
$\cal F$ and $\varphi$ from $\Lambda$.

Let $(P,P_+,[u])$ be a simple totally ordered dimension group with $P\simeq\Z^n$. Recall
that a  Bratteli diagram of $(P,P_+,[u])$ consists of
a vertex set $V$ and edge set $E$ such that $V$ is an infinite
disjoint union $V_1\sqcup V_2\sqcup\dots$, where each $V_i$ has
cardinality $n$. 
The latter condition follows from the total ordering of ${\Bbb Z}^n$.
Any pair $V_{i-1},V_i$ defines a non-empty set
$E_i\subset E$ of edges with a pair of range and source functions
$r,s$ such that $r(E_i)\subseteq V_i$ and $s(E_i)\subseteq V_{i-1}$.

An $AF$ $C^*$-algebra whose dimension group is order-isomorphic to
$(P,P_+,[u])$ is an inductive limit of multi-matrix algebras
\displaymath
\lim M_{J_1}(\C)\oplus\dots\oplus M_{J_n}(\C).
\enddisplaymath
We shall say that a Bratteli diagram $(V,E)$ {\it corresponds to}
group $(P,P_+,[u])$ if the range and source functions of $(V,E)$
represent embedding scheme of the above multi-matrix algebras. 
(In other words, an $AF$-algebra defined by $(V,E)$ has Elliott
group $(P,P_+,[u])$.)

The equivalence class of Bratteli diagrams corresponding to a
simple totally ordered 
dimension group of form $\Z^n$ has a representative $(V,E)$ with no multiple
edges, since every positive integral matrix decomposes into a finite product
of non-negative matrices whose entries are zeros and ones.
For the sake of simplicity, we always assume this
case of Bratteli diagrams.

By an {\it infinite path} on $(V,E)$ we shall mean an infinite sequence
of edges $(e_0,e_1,\dots)$ such that $e_0\in E_0,e_1\in E_1,$ etc.
The set of all infinite paths on $(V,E)$ is denoted by $X$.
Let us identify ``coordinates'' $x_i$ of $x\in X$ with vector
$(e_0,e_1,\dots)$. Fix $x,y\in X$. Metric $d(x,y)=1/2^k$, where
\displaymath
k=\max \{l\in\N~|~x_i=y_i~\hbox{for} ~i<l\}, 
\enddisplaymath
turns $X$ into an absolutely disconnected topological space
which is called a {\it Bratteli-Cantor compactum}. To construct an embedding
$X\to {\Bbb H}$ where each $x\in X$ represents a geodesic, a portion
of symbolic dynamics is needed.

\medskip\noindent
\underline{Koebe-Morse coding of geodesics}.
Let $M_{g,m}$ be a hyperbolic surface of genus $g$ with $m$ totally
geodesic boundary components $v_1,\dots,v_m$. We dissect $M_{g,m}$ to a simply
connected surface as follows \cite{Mor}. Let $P$ be
an arbitrary point of $v_m$. One draws geodesic segments $h_1,\dots,h_{m-1}$
from $P$ to some arbitrary chosen points of $v_1,\dots,v_{m-1}$. (Thus, $h_i$
have only $P$ as common point.)
Next one dissects the handles of $M_{g,m}$ by closed geodesics $c_1,\dots,c_{2g}$
issued from point $P$. Clearly, the resulting surface is simply connected
and has the boundary
\begin{equation}
c_1,\dots,c_{2g}; h_1,\dots,h_{m-1}.
\end{equation}

\bigskip
Now given geodesic half-circle $S\subset {\Bbb H}$ passing through unique
point $0\in\tau$ one relates an infinite sequence of symbols
\begin{equation}\label{code}
\sigma_1,\sigma_2,\sigma_3,\dots,
\end{equation}
which ``take values'' in the set $\sigma$. One prescribes $\sigma_p, 
~p=1,\dots,\infty$ a ``value'' $g_i,~1\le i\le n$ if and only if $S$
has a transversal intersection point with the side $a_i=b_i$ of $p$-th
image of $G_{\sigma}\in\tau$. (In other words, code (\ref{code})
``counts'' points of intersection of $S$ with ``sides'' of tessellation
$\tau$.) Sequence of symbols (\ref{code}) is called a {\it Koebe-Morse code}
of geodesic $S$.

Morse showed that there is a bijective correspondence
between sequences (\ref{code}) satisfying some admissibility requirements
\footnote{Namely, there should be no words with the syllabi $a_ib_i$ or $b_ia_i$,
where $a_i$ and $b_i$ are ``dual'' symbols from the alphabet $G_{\sigma}$.}
and the set of non-periodic geodesics on surfaces of negative curvature;
see the bibliography to Morse and Hedlund \cite{MoH}.   
\begin{lem}
Let $S$ be a geodesic with the Koebe-Morse code $(\sigma_1,\sigma_2,\dots).$
Then any congruent to $S$ geodesic $S'$ will have the same Koebe-Morse code,
except possibly in a finite number of terms.
\end{lem}
{\it Proof.} This follows from the definition of coding and invariance
of $\tau$ by the $G$-actions.
$\square$

\bigskip\noindent
Let $X$ be a Bratteli-Cantor compactum. Let $V_i\to\sigma$ be a bijection
between the vertices $V=V_1\sqcup V_2\sqcup\dots$ of $(V,E)$ and the
set of symbols $\sigma$. This bijection can be established by labeling
each element of $V_i$ from the left to the right by symbols $\{g_1,\dots,g_n\}$.
Thus, every $x\in X$ is a ``symbolic geodesic'' $(x_1,x_2,\dots)$ whose
``coordinates'' take values in $\sigma$. Each sequence is admissible and
by Morse's Theorem realized by a (class of congruent) geodesic whose
Koebe-Morse code coincides with $(x_1,x_2,\dots)$. Where there is no
confusion, we refer to $x\in X$ as a geodesic line in the complex plane
$\Bbb H$. 
\begin{lem}\label{closure}
Let $l_x$ be an image of geodesic $x\in X$ on the surface $M_{g,m}$
under projection ${\Bbb H}\to {\Bbb H}/G$. If Bratteli diagram $(V,E)$
is simple, then $l_y\in Clos~l_x$ for any $y\in X$.
\end{lem}
{\it Proof.} The simplicity of $(V,E)$ means that every infinite
path $x\in X$ is transitive, i.e. any finite ``block'' of symbols
$\{x_n,x_{n+1},\dots,x_{n+k}\}$ occurs ``infinitely many times''
in the sequence $x=(x_1,x_2,\dots)$. 
Indeed, simplicity of $(V,E)$ means the absence of non-trivial
ideals in the corresponding $AF$ $C^*$-algebra. Using Bratteli's
dictionary (\cite{Bra}) between ideals and connectedness properties of $(V,E)$,
it can be easily shown that arbitrary infinite path in $(V,E)$
``visits'' any given finite sequence of vertices infinitely often.

Suppose that $B_k$ is a block of symbols of length $k\ge1$. Let
\displaymath
x=(x_1,\dots,x_{n-1},B_k,x_{n+k+1},\dots),\qquad
y=(y_1,\dots,y_{m-1},B_k,y_{m+k+1},\dots),
\enddisplaymath
be the first time $B_k$ appears in sequences $x,y\in X$. 
By a congruent transformation, the geodesics $x,y\in {\Bbb H}$ can 
be brought to the from
\displaymath
x'=(B_k,x_{k+1},\dots),\qquad
y'=(B_k,y_{k+1},\dots).
\enddisplaymath
This means that $dist~(x',y')\le 1/2^k$. Since $B_k$ occurs in
sequences $x$ and $y$ infinitely often, lemma follows.
$\square$ 
\begin{lem}\label{lamination}
The set $Clos~l_x$ of Lemma \ref{closure}  is a set $\Lambda\subset M_{g,m}$
consisting of continuum of irrational geodesic lines.
\end{lem}
{\it Proof.} This follows from the proof of Lemma \ref{closure}.
$\square$

\bigskip\noindent
By Lemmas \ref{closure} and \ref{lamination}, $\Lambda$ is
homeomorphic to Thurston's geodesic lamination on a surface of genus
$g\ge 2$; cf Thurston \cite{Thu}. To finish the proof of theorem,
one needs ``blow-down'' $\Lambda$ to a measured foliation $\F$.
The required interval exchange transformation $\varphi$ is the
``mapping of first return'' on a global transversal to $\F$.
By the construction, $\varphi$ is minimal and has $n=2g+m-1$ intervals
of continuity.

Theorem \ref{range} is proven.
$\square$

\section{Criterion of strict ergodicity for 
the interval exchange transformations}
In general, group $(K_0,K_0^+,[1])$ may be {\it not} totally ordered.
The total order happens if and only if positive cone
$K_0^+$ is bounded by unique ``hyperplane'' in the
``space'' $K_0$. There exists up to $(n-2)$
hyperplanes in a group $K_0$ of rank $n$ which constitute
a boundary of $K_0^+$; cf Goodearl \cite{G}, p.217.

A {\it state} is a homomorphism $f$ from $(K_0,K_0^+,[1])$
to $\R$ such that $f(K_0^+)\subset\R_+$ and $f([1])=1$.
The space of states $S_{\bullet}$ is dual to the linear space $K_0^+$.
From this point of view, ``hyperplanes'' correspond to linearly
independent ``vectors'' of the space $S_{\bullet}$. A total order
is equivalent to the requirements $dim~S_{\bullet}=1$ and absence
of ``infinitesimals'', cf Effros \cite{E} p.26.

Let $\varphi$ be an interval exchange transformation built upon
$(K_0,K_0^+,[1])$. Invariant measures of $\varphi$ form a vector
space relatively sums and multiplication of measures by positive
reals. This vector space is isomorphic to $S_{\bullet}$. 
The requirement $\varphi$ to be strictly
(uniquely) ergodic is equivalent to the claim $S_{\bullet}$
be one-dimensional. Strict ergodicity of the interval exchange
transformations has been a challenging problem in the area
for years. (Find a working criterion to determine whether
given $\varphi$ is strictly ergodic.)

This saga started in 1975 when examples of interval exchange
transformations with two and three invariant ergodic measures
became known due to Keynes and Newton. Keane made an assumption
that the ``majority'' of transformations $\varphi$ are strictly
ergodic. This assumption was turned to a theorem independently
by Masur and Veech who used for this purpose the Teichm\"{u}ller
theory and topological dynamics, respectively.  
The proof of Keane's conjecture based on Theorem \ref{range} is given 
in Section 2.3.

In this section we establish strict ergodicity for a class of
interval exchange transformations which we call ``stationary''.
The name comes from  theory of ordered abelian groups, because such 
transformations have stationary Bratteli diagrams; cf Effros \cite{E}.
Foliations that correspond to such transformations are known
as {\it pseudo-Anosov} or foliations whose  leaves
 are 1-dimensional basic sets of the pseudo-Anosov homeomorphisms of a compact
surface.
\begin{dfn}
Let $\varphi=\varphi(\lambda,\pi,\varepsilon)$ be an interval exchange
transformation whose Bratteli diagram is given by the infinite
sequence of multiplicity matrices
\begin{equation}\label{sequence}
\{P_{Y_1},P_{Y_2},P_{Y_3},\dots \}.
\end{equation}
If the set (\ref{sequence})
can be divided into the blocks $B_k=\{P_{Y_n},P_{Y_{n+1}},\dots,P_{Y_{n+k}}\}$
such that $P_{Y_n}P_{Y_{n+1}}\dots P_{Y_{n+k}}=P$, then $\varphi$ is called stationary.
In particular, $\varphi$ is stationary if $P_{Y_1}=P_{Y_2}=P_{Y_3}=\dots=P$.
\end{dfn}
\begin{thm}\label{strict}
Every stationary interval exchange transformation $\varphi$ is strictly
ergodic.
\end{thm}
{\it Proof.}  The proof is based on the 
Perron-Frobenius Theorem. A dual (projective) limit
\begin{equation}\label{eq21}
(\R^n)^*\buildrel\rm P_{Y_1}\over\longrightarrow (\R^n)^*
\buildrel\rm P_{Y_2}\over\longrightarrow
(\R^n)^*\buildrel\rm P_{Y_3}\over\longrightarrow \dots
\end{equation}
consists of operators $P_{Y_i}$ acting
on the dual space $(\R^n)^*$ to $\Z^n\subset\R^n$. (In other words,
we identify the space of positive homomorphisms $\Z^n\to\R$ and 
the space of linear functionals $\R^n\to\R$.) The diagram (\ref{eq21})
converges to the state space $S_{\bullet}$ of dimension group
$(K_0,K_0^+,[1])$. When $P_{Y_i}=P$, where $P$ is a matrix with
strictly positive entries, or can be reduced to this case, then
there exists a maximal simple eigenvalue $\lambda>0$ of matrix $P$
(Perron-Frobenius Theorem). The eigenvector $x_{\lambda}$ defines
a 1-dimensional $P$-invariant subspace of $(\R^n)^*$ lying in the
limit of diagram (\ref{eq21}) and which is identified with $S_{\bullet}$.
Let us formalize this idea.

For $1\le i\le n$ denote by $e_i$
and $e_i^*$ the vectors of canonical bases in the vector space 
$\R^n$ and the dual space $(\R^n)^*$. 
By $(\alpha_1,\dots,\alpha_n)$ and $(\alpha_1^*,\dots,\alpha_n^*)$
we denote vectors in $\R^n$ and $(\R^n)^*$. 
$(\R^n)^+$ and $(\R^n)^{*+}$ are collections of vectors
whose coordinates are $\alpha_i\ge0$ and $\alpha_i^*\ge0$, respectively.
The same notation $(\Z^n)^+$ is reserved for the integer vectors of $\R^n$.  
By $\Delta_0\subset (\R^n)^*$ we understand n-dimensional simplex
spanned by the vectors $0, e_1^*,\dots,e_n^*$. 
To each linear mapping $\phi: \R^n\to\R^n$ one associates
a dual mapping $\phi^*:(\R^n)^*\to (\R^n)^*$.

Denote by $P$ an dimension group corresponding to the limit
\begin{equation}\label{eq24}
P\cong\lim_{k\to\infty}P_k,
\end{equation}
where $P_k$ are ordered groups whose positive cone is defined
to be an inverse of $k$-th iteration of the set $(\Z^n)^+$
under the automorphism $\phi$:
\begin{equation}\label{eq25}
P_k^+=\phi^{-k}[(\Z^n)^+].
\end{equation}
The set $\Delta_0$ has been introduced earlier. For $k=1,\dots,\infty$
we let
\begin{equation}\label{eq26}
\Delta_k=S_{\bullet}(P_k),
\end{equation}
which is a state space of the group $P_k$.
Define $\Delta_k$ to be a simplex spanned by the vectors 
$0,J_1(k),\dots,J_n(k)$, where 
\begin{equation}\label{eq27}
J_1(k)={X_1(k)\over ||X_1(k)||},\quad\dots\quad, J_n(k)={X_n(k)\over ||X_n(k)||},
\end{equation}
and 
\begin{equation}\label{eq28}
X_1(k)=\phi^k(e_1^*),\quad\dots\quad, X_n(k)=\phi^k(e_n^*).
\end{equation}
It is evident that
\begin{equation}\label{eq29}
\Delta_0\supseteq\Delta_1\supseteq\Delta_2\supseteq\dots\supseteq\Delta_{\infty},
\end{equation}
where 
\begin{equation}\label{eq30}
\Delta_{\infty}=\bigcap_{k=1}^{\infty}\Delta_k.
\end{equation}
A recurrent formula linking $X_i(k-1)$ and 
$X_i(k)$ is given by the equation:
\begin{equation}\label{eq31}
X_i(k)=\sum_{j=1}^np_{ij}X_j(k-1),
\end{equation}
where $p_{ij}$ are the entries of matrix of ``partial multiplicities'' $P_y$.

Note that simplex $\Delta_{\infty}$ has
dimension $r\le n$.
The original problem of calculating the state space $S_{\bullet}$
is reduced  to
calculation of the asymptotic simplex $\Delta_{\infty}$
whose spanning vectors are linked by equation (\ref{eq31}).
We shall see that $\Delta_{\infty}$ can be completely calculated under
hypothesis of Theorem \ref{strict}. The following lemma is basic.
\begin{lem} {\bf (Perron-Frobenius)}\label{lm6}
A strictly positive $n\times n$ matrix $P=(p_{ij})$ always has a real
and positive eigenvalue $\lambda$ which is a simple root of the characteristic
equation and exceeds the moduli of all the other characteristic values.
To this maximal eigenvalue $\lambda$ there corresponds an eigenvector
$x_{\lambda}=(x_{\lambda}^1,\dots,x_{\lambda}^n)$ with positive 
coordinates $x_{\lambda}^i>0, ~i=1,\dots,n$.
\footnote{In fact, there exists more general statement due to Frobenius
which treats matrices with non-negative entries. Because of exceptional
importance of this statement in understanding why the unique ergodicity
may vanish, and also due to the clear connection of Frobenius theorem
with the root systems of Coxeter-Dynkin, we give the formulation of this
theorem in below.

\noindent
{\bf Theorem (Frobenius)} {\it ~An irreducible non-negative $n\times n$
matrix $P=(p_{ij})$ always has a positive eigenvalue $\lambda$ that
is a simple root of the characteristic equation. The moduli of all
the other eigenvalues do not exceed $\lambda$. To the maximal eigenvalue
$\lambda$ there corresponds an eigenvector with positive coordinates.

Moreover, if $P$ has $r$ eigenvalues $\lambda_0=\lambda,~\lambda_1,
\dots,\lambda_{r-1}$ of modulus $\lambda$, then these numbers are
all distinct and are roots of the equation $$z^r-\lambda^r=0.$$
More generally: The whole spectrum $\lambda_0,\lambda_1,\dots,\lambda_{n-1}$
of $P$, regarded as a system of points in the complex plane, goes over
into itself under a rotation of the plane by the angle $2\pi/r$.
If $r>1$ then $P$ can be put into the following cyclic normal form 
$$\left(\matrix{O & P_{12} & O      & \dots & O\cr
                O & O      & P_{23} & \dots & O\cr
           \vdots &        &        &       & \vdots \cr
                O & O      & O      & \dots & P_{r-1,r}\cr
           P_{r1} & O      & O      & \dots & O }  
\right)$$ 
where $P_{ij}$ are non-zero square blocks along the main diagonal
and $O$ are zero square blocks elsewhere.}}
\end{lem}
{\it Proof of lemma.} Let $x=(x_1,x_2,\dots,x_n)$ be a fixed
vector. A function 
\begin{equation}\label{eq32}
r_x=\min_{1\le i\le n}{(Px)_i\over x_i}
\end{equation}
is introduced. We have $r_x\ge 0$ since 
\begin{equation}\label{eq33}
(Px)_i=\sum_{j=1}^np_{ij}x_j,
\end{equation}
is a non-negative matrix. However, in the definition of minimum (\ref{eq32}),
the values of $i$ for which $x_i=0$ are excluded. 
The lemma follows from the variational principle for the
continuous function $r_x$, which must assume a maximal
value for some vector $x$ with non-negative coordinates.
$\square$

\bigskip\noindent
Now we can finish the proof of main theorem. By strict positivity
of matrix $P$ and Lemma \ref{lm6}, there is a positive maximal eigenvalue
$\lambda$, whose eigenvector $x_{\lambda}$ has positive coordinates. 
Notice that 
\begin{equation}\label{eq34}
\phi^k(x_{\lambda})=(\lambda)^kx_{\lambda},
\end{equation}
so that the iterations of $\phi$ leave invariant a 1-dimensional
linear subspace $\{\alpha\}$ spanned by $x_{\lambda}$.
All other vectors in $(\R^n)^{*+}$ converge accordingly 
(\ref{eq31}) to the subspace $\{\alpha\}$.
We conclude that 
\begin{equation}\label{eq35}
\Delta_{\infty}=\{\alpha\},
\end{equation}
which is  one-dimensional. Theorem \ref{strict} is proved. 
$\square$     
\begin{rmk}
In the context of ordered abelian groups, Theorem \ref{strict} 
was known to Effros \cite{E} and Elliott \cite{Ell}. 
The main ingredients of proof can be traced in the work \cite{Vee1} of
W.~Veech.
\end{rmk}

\section{Masur-Veech Theorem}
Theorem of Masur and Veech is formulated in Section 5.2. There are
two known proofs of this theorem, due to Masur \cite{Mas} who
used complex analysis and Teichm\"{u}ller theory and Veech \cite{Vee}
who used methods of topological dynamics. In this section we
suggest an independent proof using Theorem \ref{range} and
one lemma of Morse and Hedlund from symbolic dynamics; cf Morse
and Hedlund \cite{MoH}. 

\bigskip\noindent
\underline{Parametrization of $(K_0,K_0^+,[1])$}.
\footnote{The idea of such parametrization was communicated
to the author by G.~A.~Elliott.}
Let $\Bbb H,\tau$ and $S$ be as in Definition \ref{positivecone}
of Section 2. Without loss of generality we assume that $S$ is
a unit semi-circle in the complex plane $\Bbb H$. Consider a family
$S_t$ of the unit semi-circles parametrized by real numbers equal
to a ``horizontal shift'' of $S$ in $\Bbb H$. (In other words,
$t$ is equal to the $x$-coordinate of the centre of unit circle $S_t$.)
A family of dimension groups which are defined by ``positive cones''
$S_t$, we shall denote by $(P,P_+^t,[u])$. By results of Sections 1-3
every dimension group of form $\Z^n$ has a representative in $(P,P_+^t,[u])$
and every measured foliation (with fixed singularity data) arises
in this way.
\begin{thm}\label{Masur-Veech}
Denote by $\F_t$ a family of measured foliations corresponding to
$(P,P_+^t,[u])$ and by $t_1\sim t_2$ an equivalence relation on $\R$
identifying topologically equivalent foliations $\F_{t_1}$ and $\F_{t_2}$.
If $X=\R/\sim$ is a topological space, then for a residual set of
the second category in $X$ foliation $\F_t$ is strictly ergodic.  
\end{thm}
{\it Proof.} The idea is to apply Koebe-Morse coding to 
each geodesic $S\in S_t$. In this setting, $X$ becomes a space
of symbolic sequences with the topology described in Section 2.
It is not hard to see that strict ergodicity of individual geodesic $S$
 is equivalent to ``uniform approximation'' of $S$ by periodic
 sequences of length $N$. (Using the terminology of Morse and
 Hedlund, such approximation property of a geodesic means that
 a transitivity index $\phi(N)$ tends to a  covering
 index $\theta(N)$ of the geodesic as $N\to\infty$; cf Morse and
 Hedlund \cite{MoH}.) The same authors proved that 
 $\lim_{\N\to\infty}\inf~{\phi(N)\over\theta(N)}=1$ for a residual
 set of the second category in $X$. Let us give the details of
 this construction.

Let $S$ be a geodesic in the complex plane $\Bbb H$ and
\displaymath
\sigma_1,\sigma_2,\sigma_3,\dots,
\enddisplaymath
the Koebe-Morse code of $S$ which we shall call a {\it ray};
cf Section 3. The ray $R$ is {\it transitive} if it contains
a copy of each admissible block. (Block is shorthand for a finite
sequence of symbols.) Function $\phi:\N\to\N$
of a transitive ray $R$ is called a {\it transitivity index}
if the initial block of $R$ of the length $\phi(N)$ contains
all admissible blocks of the length $N$ and there are no shorter
initial subblocks with this property. Dropping the claim that
block $B\subset R$ is initial gives us function $\theta:\N\to\N$ 
which is called a {\it covering index} of the recurrent ray $R$.
These functions satisfy an obvious inequality:
\displaymath
\phi(N)\ge\theta(N).
\enddisplaymath
\begin{lem} {\bf (Morse-Hedlund)}\label{Morse-Hedlund}
A set of rays whose transitivity index and covering index
satisfy the condition
\displaymath
\lim_{N\to\infty}\inf {\phi(N)\over\theta(N)}=1,
\enddisplaymath
is a residual set of the second category in the space $X$ of
all infinite rays endowed with the topology described in Section 2. 
\end{lem}
{\it Proof of lemma.} For a complete proof see \cite{MoH}.
Denote by $Y$ the set of rays satisfying the condition of lemma.
The following items will be proved consequently:

\medskip
(i) $Y$ is not empty;

\smallskip
(ii) $Y$ is everywhere dense in $X$;

\smallskip
(iii) The complement of $Y$ is nowhere dense in $X$.

\bigskip
(i) Let $H(n)$ be a block of minimum length containing all admissible
blocks of length $n$. For a growing sequence of integers $r_0,r_1,\dots,r_{k-1}$
consider a block
\begin{equation}\label{block}
H(r_0)\sigma_1 H(r_1)\sigma_2\dots\sigma_{k-1} H(r_{k-1})\sigma_k
\end{equation}
of length $m_k$. Let us choose $r_k$ sufficiently large so that
\displaymath
{\theta(r_k)+m_k\over\theta(r_k)}<1+\delta_k,
\enddisplaymath
where $\delta_k$ is a vanishing positive real. The transitivity index of (\ref{block})
satisfies the inequality
\displaymath
\phi(r_k)\le \theta(r_k)+m_k.
\enddisplaymath
By the construction, ${m_k\over\theta(r_k)}\to 0$ as $k\to\infty$ so that (\ref{block})
satisfies the condition of lemma.

\medskip
(ii) Let $A$ be an arbitrary admissible block of length $k$ and $R\in Y$. 
For a suitably chosen $\sigma$ the ray
\displaymath
R'=A\sigma R
\enddisplaymath
is admissible. We have the following inequalities for $R$ and $R'$:
\displaymath
\theta(N)\le \phi'(N)\le\phi(N)+k+1,
\enddisplaymath
where $\phi$ and $\phi'$ are the transitivity indices of $R$ and $R'$.
The condition of lemma is satisfied and therefore $R'\in Y$.  

\medskip
(iii) This item follows from (ii) and accurate construction of closed
sets lying in the complement of $Y$; cf \cite{MoH} for the details.
This argument finishes the proof of Morse-Hedlund lemma.
$\square$

\bigskip\noindent
Let $R$ be a transitive ray and $B_1,\dots,B_k$ admissible blocks
of length $N$. We say that $R$ is {\it uniformly distributed} relatively
$B_1,\dots,B_k$ if
\displaymath
\phi(N)=kN.
\enddisplaymath
(In other words, each admissible block appears in the initial block of $R$
with the ``probability'' $1/k$.) In the geometric terms this means that geodesic
$R$ is located at the same distance from  periodic geodesics
\displaymath
B_1,B_1,\dots;\quad B_2,B_2,\dots;\quad\dots;\quad B_k,B_k,\dots.
\enddisplaymath
\begin{lem}\label{dist}
Suppose that $R$ is uniformly distributed relatively admissible blocks
$B_1,\dots,B_k$ for each integer $N>0$. Then
\displaymath
\lim_{N\to\infty}\inf {\phi(N)\over\theta(N)}=1.
\enddisplaymath
\end{lem}
{\it Proof of lemma.} This follows from the equality $\phi(N)=\theta(N)$.
$\square$

\bigskip\noindent
To finish the proof of Theorem \ref{Masur-Veech} it remains to notice that
strict ergodicity of $R$ is equivalent to uniform distribution of periodic
``blocks'' in $R$ and apply Lemmas \ref{Morse-Hedlund} and \ref{dist}.
$\square$

\section{Appendix}
\subsection{Interval exchange transformations}
Let $n\ge2$ be a positive integer and let $\lambda=(\lambda_1,\dots,\lambda_n)$
be a vector with positive components $\lambda_i$ such that 
$\lambda_1+\dots+\lambda_n=1$. One sets 
\displaymath    
\beta_0=0,\qquad \beta_i=\sum_{j=1}^i\lambda_j,\qquad    
v_i=[\beta_{i-1},\beta_i)\subset [0,1].    
\enddisplaymath  
Let $\pi$ be a permutation on the index set $N=\{1,\dots,n\}$
and $\varepsilon=(\varepsilon_1,...,\varepsilon_n)$ a vector
with coordinates $\varepsilon_i=\pm 1, i\in N$.
An {\it interval exchange transformation} is a mapping 
$\varphi(\lambda,\pi,\varepsilon): [0,1]\to [0,1]$
which acts by piecewise isometries  
\displaymath    
\varphi (x)=\varepsilon_i x-\beta_{i-1}+\beta^{\pi}_{\pi(i)-1}, 
\qquad x\in v_i,    
\enddisplaymath
where $\beta^{\pi}$ is a vector corresponding to
$\lambda^{\pi}=(\lambda_{\pi^{-1}(1)},\lambda_{\pi^{-1}(2)},...,    
\lambda_{\pi^{-1}(n)})$. Mapping $\varphi$ preserves or reverses
orientation of $v_i$ depending on the sign of $\varepsilon_i$.
If $\varepsilon_i=1$ for all $i\in N$ then the interval exchange    
transformation is called {\it oriented}. Otherwise, the interval exchange 
transformation is said to have {\it flips}.

Interval exchange transformation is said to be {\it irreducible}
if $\pi$ is an irreducible permutation. An irreducible interval exchange 
transformation $T$ is called {\it irrational} if the only rational relation 
between numbers $\lambda_1,...,\lambda_n$ is given by the equality $\lambda_1+...+    
\lambda_n=1$. Recall that measure $\mu$ on $[0,1]$ is called invariant,
if $\mu(\varphi(A))=\mu(A)$ for any measurable subset $A\subseteq [0,1]$. 
The following theorem due to M.~Keane (\cite{Kea}) estimates the number of 
invariant measures.

\medskip\noindent
{\bf Finiteness Theorem.}
{\it Let $\varphi$ be an irrational interval    
exchange transformation of $n$ intervals. Then there are at most    
finitely many ergodic invariant measures 
whose number cannot exceed $n+2$.    
}

\bigskip\noindent    
In case of the interval exchange transformations without flips, there exists
an estimate of the number of invariant ergodic measures due to Veech \cite{Vee}.

\bigskip\noindent
{\bf Veech Theorem} 
{\it Let $\varphi$ be an irrational interval exchange transformation without flips 
on $n\ge 2$ intevals. Then the number of invariant ergodic measures of $\varphi$
is less or equal to $[{n\over 2}]$, where $[\bullet]$ is integer part of the number.
}

\subsection{Measured foliations}
Measured foliations are suspensions over the interval    
 exchange transformations which preserve the ergodic measure on    
 intervals and such that their singularity set consists of $p$-prong
 saddles, $p\ge 3$. Measured foliations can be    
 defined via closed 1-forms which is more elegant way     
 due to Hubbard, Masur and Thurston.    
 \begin{dfn}    
 {\bf (Hubbard-Masur-Thurston)} Let $M$ be a compact $C^{\infty}$    
 surface of genus $g>1$, without boundary. A measured foliation    
 $\F$ on $M$ with singularities of order $k_1,...,k_n$ at    
 points $x_1,...,x_n$ is given by an open cover    
 $U_i$ of $M\backslash\{x_1,...,x_n\}$ and non-vanishing    
 $C^{\infty}$ real valued closed 1-form $\phi_i$ on each    
 $U_i$, such that    
    
 (i) $\phi_i=\pm \phi_j$ on $U_i\cap U_j$;    
    
 (ii) at each $x_i$ there is a local chart $(u,v):V\rightarrow\R^2$    
 such that for $z=u+iv$, $\phi_i=$ Im $(z^{k_i/2}dz)$ on    
 $V\cap U_i$, for some branch of $z^{k_i/2}$ in $U_i\cap V$.    
    
 Pairs $(U_i,\phi_i)$ are called an atlas for $\F$.    
 \end{dfn}
 As it follows from the definition, apart from the singular points,    
 measured foliations look like a non singular volume preserving    
flows. In singularities, the substitution $z\mapsto re^{i\psi}$    
brings $\phi_i$, mentioned in (ii), to the form    
\begin{displaymath}    
\phi_i=r^{{k_i\over 2}}[\sin({k_i\over 2}+1)\psi dr+r\cos    
({k_i\over 2}+1)\psi d\psi].    
\end{displaymath} 
It can be readily established, that $\phi_i$ are closed differential    
1-forms, that is $d\phi_i=0$ for all $k_i\ge 1$.    
To see what  singularities are generated by the above formula,    
let us consider a vector field $v_i$, given by the system of differential    
equations    
\begin{displaymath}    
{dr\over dt} = -r\cos({k_i\over 2}+1)\psi,\qquad    
{d\psi\over dt} = \sin ({k_i\over 2}+1)\psi.    
\end{displaymath}
Clearly, $v_i$ is tangent to a foliation given by the equation    
$\phi_i=0$. Our prior interest is to study the behavior of trajectories    
of $v_i$ in a narrow stripe $\Pi=\{(r,\psi)|-\varepsilon\le r\le \varepsilon,    
0\le\psi\le2\pi\}$. There are exactly $k_i+2$ equilibria $p_n\in\Pi$,    
which have the coordinates $(0,{2\pi n\over k_i+2})$, where $n\in\N$    
varies from $0$ to $k_i+2$. The linerization of the vector field $v_i$    
in these points yields    
\begin{displaymath}    
  A(p_n)=\left(\matrix{(-1)^{n+1}&0\cr    
                   & \cr    
                 0&(-1)^n({k_i\over 2}+1)\cr}\right)    
\end{displaymath}
Therefore all $p_n$ are the saddle points.    
One maps the half-stripe $r\ge0$ to the neighbourhood of the singular    
point $x_i$. Generally, a singular    
point $x_i$ of the order $k_i$ is a $(k_i+2)$-prong saddle of a measured    
foliation $\F$.

\bigskip
Let $M$ be a compact surface and $\F$ a measured foliation on    
$M$. By {\it measure} $\mu$ of $\F$ one understands a line element    
$||\phi||$ related with the point $x\in M$, induced in each $x\in U_i$    
by $||\phi_i(x)||$. It measures a `transversal length' of $\F$, since    
$\mu$ vanishes in direction tangent to the leaves of $\F$.

Take  a cross-section to the 
measured foliation $\F$. $\F$ induces an interval exchange transformation 
$\varphi$  on this cross-section. Depending on orientability of $\F$,
$\varphi$ may have flips.  Flips are excluded if $\F$ is an orientable
measured foliation (in this case $\F$ is given by orbits of a measure-preserving
flow). For orientable measured foliations, an estimate of number of invariant 
ergodic measures is due to Sataev \cite{Sat}.

\bigskip\noindent
{\bf Sataev Theorem}
{\it Let $n$ and $k$ be a pair of natural numbers, such that    
$n\ge k$ and let $M$ be a compact orientable surface of genus $n$.    
Then there exists a $C^{\infty}$ orientable measured foliation    
$\F$ on $M$ whose singularity set consists of 4-separatrix saddles and 
 which has exactly $k$ invariant ergodic measures.    
}

\bigskip\noindent
An important question arises when measured foliation has 
a unique invariant measure. It was conjectured by M.~Keane and proved 
by H.~Masur and W.~Veech that `almost all'    
measured foliations have a unique    
invariant measure, which is a multiple of  Lebesgue measure.

\bigskip\noindent
{\bf Masur-Veech Theorem} 
{\it (\cite{Mas},\cite{Vee})
Suppose that family $\F_t$ of measured foliations 
is given by trajectories of a holomorphic quadratic differential 
$e^{it}\phi$ on the surface $M$. Then for `almost all' values of $t$ foliation
$\F_t$ is strictly ergodic.
}

\subsection{${\cal O}_{\lambda}$ as a crossed product $C^*$-algebra}
\begin{lem}\label{groups}
Let $\varphi=\varphi(\lambda,\pi,\varepsilon)$ be an interval
exchange transformation and $\lambda=(\lambda_1,\dots,\lambda_n)$.
Then $K_0({\cal O}_{\lambda})=\Z^n$ and $K_1({\cal O}_{\lambda})=\Z$.
\end{lem}
{\it Proof.} 
Let $p_1,\dots, p_n$ be the set of discontinuous points of the mapping
$\varphi$. Denote by $Orb~\varphi=\{\varphi^m(p_i):~1\le i\le n,~n\in\Z \}$ a set
of full orbits of these points. When $\varphi$ is irrational, the set 
$Orb~\varphi$ is a dense subset in $[0,1]$. 
We  replace every point $x\in Orb~\varphi$ in the interior of $[0,1]$ by
two points $x^-<x^+$ moving apart banks of the cut. 
The obtained set is a Cantor set denoted by $X$.

A mapping $\varphi: X\to X$ is defined to coincide with
the initial interval exchange transformation on  
$[0,1]\backslash Orb~\varphi\subset X$
 prolonged to a homeomorphism of $X$. Mapping $\varphi$
is a minimal homeomorphism of $X$, since there
are no proper, closed, $\varphi$-invariant subsets of $X$ except
the empty set. Thus, ${\cal O}_{\lambda}=C(X)\rtimes_{\varphi}\Z$
is a crossed product $C^*$-algebra, where $C(X)$ denotes a $C^*$-algebra
of continuous complex-valued functions on $X$. The following diagram of
Pimsner and Voiculescu consists of exact sequences:

\bigskip
\begin{picture}(400,100)(0,0)
\put(85,83){\vector(3,0){50}}
\put(30,80){$K_0(C(X))$}
\put(8,20){$K_1(C(X)\rtimes_{\varphi}\Z)$}
\put(140,80){$K_0(C(X))$}
\put(195,83){\vector(3,0){40}}
\put(240,80){$K_0(C(X)\rtimes_{\varphi}\Z)$}
\put(90,90){$id-\varphi_*$}
\put(215,90){$i_*$}
\put(55,35){\vector(0,3){35}}
\put(275,72){\vector(0,-3){35}}
\put(140,20){$K_1(C(X))$}
\put(135,23){\vector(-3,0){40}}
\put(240,23){\vector(-3,0){40}}
\put(250,20){$K_1(C(X))$}
\put(115,30){$i_*$}
\put(210,30){$id-\varphi_*$}
\end{picture}

\bigskip\noindent
It was proved in \cite{Put1} that $K_0(C(X))\simeq\Z^n$
and $K_1(C(X))\simeq 0$. To obtain the conclusion of Lemma \ref{groups} it remains to 
calculate all short exact sequences in the diagram of Pimsner and
Voiculescu.
$\square$

\subsection{Rotation numbers}
One of the striking invariants of algebra ${\cal O}_{\lambda}$
are rotation numbers associated to this algebra. In the dynamical
context, rotation numbers are equal to ``average inclination''
of leaves of measured foliation relatively a coordinate system
on $M_n$. (In fact, the original study of ${\cal O}_{\lambda}$
was motivated by possibility to introduce such numbers; cf \cite{Nik2}.)
Rotation numbers for ${\cal O}_{\lambda}$ play the same role as
real numbers $\theta$ for the irrational rotation algebra $A_{\theta}$.

Recall that the cone $K_0^+\subset {\Bbb H}$ is a limit of ``rational''
cones $P_k^+\subset {\Bbb H}$:
\displaymath
K_0^+=\lim_{k\to\infty}P^+_k.
\enddisplaymath
Each $P_k^+$ is represented by a periodic geodesic $\gamma_k$.
Suppose that $g_k\in G$ is an isometry which moves geodesic $\gamma_{k-1}$
to geodesic $\gamma_k$ and let 
\displaymath
g_k=\left(\matrix{ a_i & b_i\cr c_i   & d_i}\right)\in
PSL_2(\Z)
\enddisplaymath
be an integral matrix with non-negative entries and determinant $\pm 1$, corresponding to $g_k$.
The continued fraction

\bigskip
\displaymath
\theta_{\lambda}=
{a_1\over c_1}-{c_1^{-2}\over\displaystyle {d_1\over c_1}+{a_2\over c_2}-
{\strut c_2^{-2}\over\displaystyle {d_2\over c_2}+{a_3\over c_3}\displaystyle -\dots}}
\enddisplaymath

\bigskip\noindent
converges to a real number $\theta_{\lambda}$ which is called a
{\it rotation number} associated to algebra ${\cal O}_{\lambda}$.
The importance of rotation numbers is stipulated by the following theorem.

\bigskip\noindent
{\bf Theorem}
{\it Let ${\cal O}_{\lambda}$ and  ${\cal O}_{\lambda}'$
be two $C^*$-algebras whose rotation numbers are $\theta_{\lambda}$
and $\theta_{\lambda}'$. Then ${\cal O}_{\lambda}$ is Morita equivalent to
${\cal O}_{\lambda}'$  if and only if $\theta_{\lambda}$
and $\theta_{\lambda}'$ are modular equivalent:
\displaymath
\theta_{\lambda}'={a\theta_{\lambda}+b\over c\theta_{\lambda}+d},
\qquad a,b,c,d\in\Z,\qquad ad-bc=\pm 1.
\enddisplaymath
}

\bigskip\noindent
{\it Proof.} This was proved in \cite{Nik2}. 
$\square$
\begin{cor}\label{surd}
Suppose that $\varphi$ is a stationary interval exchange
transformation described in Section 3. Then the rotation
number $\theta_{\lambda}$ is a quadratic surd (i.e. irrational
root of a quadratic equation). 
\end{cor}
{\it Proof.} By the results of Section 3, the dimension group
of ${\cal O}_{\lambda}$ is stationary and must correspond
to a periodic continued fraction (i.e. $g_1=g_2=\dots=Const$).
These fractions generate a field of quadratic algebraic numbers.
$\square$

\section*{Conclusions and open problems}
The criterion of
strict ergodicity of Section 3 is highly constructive and can be used in practice
to check whether given interval exchange transformation is strictly
ergodic or not. (This can find applications in the theory of billiards 
in the rational polygons.) Of course, these conditions are only sufficient. 
The necessary conditions seem to be an open problem so far.

Another open problem is to relate the arithmetic of rotation numbers 
$\theta_{\lambda}$ with the number of invariant measures of the transformation
$\varphi$. (In the case of strictly ergodic $\varphi$ the answer is given
by Corollary \ref{surd}.)

    
\end{document}